# SELFADJOINT EXTENSIONS OF A SINGULAR DIFFERENTIAL OPERATOR

## E. BAIRAMOV[1)], R. Ö. MERT[2)] & Z. ISMAILOV[2)]

[1)] Ankara University, Faculty of Sciences, Department of Mathematics, 06100, Ankara, TURKEY

[2)] Karadeniz Technical University, Faculty of Sciences, Department of Mathematics, 61080, Trabzon, TURKEY

e-mail: bairamov@science.ankara.edu.tr; rukiye-55@hotmail.com; zameddin@yahoo.com

### Abstract

In this work, firstly in the Hilbert space of vector-functions $L^2(H, (-\infty, a) \cup (b, +\infty)), a < b$ all selfadjoint extensions of the minimal operator generated by linear singular symmetric differential expression $l(\cdot) = i\frac{d}{dt} + A$ with a selfadjoint operator coefficient $A$ in any Hilbert space $H$, are described in terms of boundary values. Later structure of the spectrum of these extensions is investigated.

## 1. Introduction

Many problems arising in the modelling of processes of multi-particle quantum mechanics, quantum field theory, in the physics of rigid bodies and ets support to study selfadjoint extension of symmetric differential operators in direct sum of Hilbert spaces([1-3]).

The general theory of selfadjoint extensions of symmetric operators in any Hilbert space and their spectral theory tall and well-built have been investigated by many mathematicians ( [4-7] ). Applications of this theory to two point differential operators in Hilbert space of functions are continied today even.

It is known that for the existence of selfadjoint extension of the any linear closed densely defined symmetric operator $B$ in a Hilbert space $\mathcal{H}$, the necessary and sufficient condition is a equality of deficiency indices $m(B) = n(B)$, where $m(B) = \dim \ker (B^* + i)$, $n(B) = \dim \ker (B^* - i)$.

The table is changed in the multipoint case in the following sense. Let $B_1$ and $B_2$ be minimal operators generated by the linear differential expression $i\frac{d}{dt}$ in the Hilbert space of functions $L^2(-\infty, a)$ and, $L^2(b, +\infty), a < b$, respectively. In this case it is known that

$$\begin{aligned}(m(B_1), n(B_1)) &= (0,1), \\ (m(B_2), n(B_2)) &= (1,0).\end{aligned}$$

Consequently, $B_1$ and $B_2$ are maximal symmetric operators, but are not a selfadjoint.

However, direct sum $B = B_1 \oplus B_2$ of operators in a direct sum $\mathcal{H} = L^2(-\infty, a) \oplus L^2(b, +\infty)$ spaces have an equal defect numbers (1,1). Then by the general theory it has a selfadjoint extension. On the other hand it can be easily shown in the form that
$$u_2(b) = e^{i\varphi} u_1(a), \varphi \in [0, 2\pi), u = (u_1, u_2), u_1 \in D(B_1^*), u_2 \in D(B_2^*).$$

Note that a space of boundary values has an important role in the theory of selfadjoint extensions of the linear symmetric differential operators ( [6,7] ).

Let $B: D(B) \subset \mathcal{H} \to \mathcal{H}$ be a closed densely defined symmetric operator in the Hilbert space $\mathcal{H}$, having equal finite or infinite deficiency indices. A triplet $(\mathfrak{H}, \gamma_1, \gamma_2)$, where $\mathfrak{H}$ is a Hilbert space, $\gamma_1$ and $\gamma_2$ are linear mappings of $D(B^*)$ into $\mathfrak{H}$, is called a space of boundary values for the operator $B$ if for any $f, g \in D(B^*)$
$$(B^*f, g)_{\mathcal{H}} - (f, B^*g)_{\mathcal{H}} = (\gamma_1(f), \gamma_2(g))_{\mathfrak{H}} - (\gamma_2(f), \gamma_1(g))_{\mathfrak{H}},$$
while for any $F_1, F_2 \in \mathfrak{H}$, there exists an element $f \in D(B^*)$, such that $\gamma_1(f) = F_1$ and $\gamma_2(f) = F_2$.

Note that any symmetric operator with equal deficiency indices have at least one space of boundary values ([6]).

In this work in second section all selfadjoint extensions of the minimal operator generated by multipoint symmetric differential operator of first order in the direct sum of Hilbert spaces $L^2(H, (-\infty, a)) \oplus L^2(H, (b, +\infty)), a < b$ in terms of boundary values are described.

In third section the spectrum of such extensions is researched.

## 2. Description of Selfadjoint Extensions

Let $H$ be a separable Hilbert space and $a, b \in \mathbb{R}, a < b$. In the Hilbert space $L^2(H, (-\infty, a)) \oplus L^2(H, (b, +\infty))$ of vector-functions considers the following linear multipoint differential expression in form
$$l(u) = (l_1(u_1), l_2(u_2)) = (iu_1' + Au_1, iu_2' + Au_2), \ u = (u_1, u_2),$$
where $A: D(A) \subset H \to H$ is linear selfadjoint operator in $H$. In the linear manifold $D(A) \subset H$ introduce the inner product in form
$$(f, g)_+ := (Af, Ag)_H + (f, g)_H, \ f, g \in D(A).$$
Then $D(A)$ is a Hilbert space under the positive norm $\|\cdot\|_+$ respect to Hilbert space $H$. It is denoted by $H_+$. Denote the $H_-$ a Hilbert space with negative norm. It is clear that a operator $A$ is continuous from $H_+$ to $H$ and it's adjoint operator $\tilde{A}: H \to H_-$ is a extension of the operator $A$. On the other hand, the operator $\tilde{A}: D(A) = H \subset H_{-1} \to H_{-1}$ is a linear selfadjoint.

From this define by
$$\tilde{l}(u) = (\tilde{l}_1(u_1), \tilde{l}_2(u_2)), \tag{2.1}$$
where $u = (u_1, u_2)$ and $\tilde{l}_1(u_1) = iu_1' + \tilde{A}u_1, \tilde{l}_2(u_2) = iu_2' + \tilde{A}u_2$.

The minimal $L_{10}(L_{20})$ and maximal $L_1(L_2)$ operators generated by differential expression $\tilde{l}_1(\tilde{l}_2)$ in $L^2(H, (-\infty, a))(L^2(H, (b, +\infty)))$ have been investigated in [8].

The operators defined by $L_0 = L_{10} \oplus L_{20}$ and $L = L_1 \oplus L_2$ in the space $L^2 = L^2(H, (-\infty, a)) \oplus L^2(H, (b, +\infty))$ are called minimal and maximal (multipoint) operators generated by the differential expression (2.1), respectively. Note that the operator $L_0$ is a symmetric and $L_0^* = L$ in $L^2$. On the other hand, it is clear that,
$$m(L_{10}) = 0, \ n(L_{10}) = \dim H,$$
$$m(L_{20}) = \dim H, \ n(L_{20}) = 0.$$

Consequently, $m(L_0) = n(L_0) > 0$. So the minimal operator $L_0$ has a selfadjoint extension ([4]). For example, the differential expression $l(u)$ with a boundary condition $u(a) = u(b)$ generates a selfadjoint operator in $L^2$.

Here it is described all selfadjoint extensions of the minimal operator $L_0$ in $L^2$ in terms of the boundary values.

In first note that the following proposition which validity of this cleam can be easily proved.

**Theorem 2.1:** The triplet $(H, \gamma_1, \gamma_2)$,
$$\gamma_1: D(L_0^*) \to H, \quad \gamma_1(u) = \frac{1}{i\sqrt{2}}(u_1(a) + u_2(b)),$$
$$\gamma_2: D(L_0^*) \to H, \quad \gamma_2(u) = \frac{1}{\sqrt{2}}(u_1(a) - u_2(b)), \quad u = (u_1, u_2) \in D(L_0^*)$$
is a space of boundary values of the minimal operator $L_0$ in $L^2$.

**Proof:** For the arbitrary $u = (u_1, u_2)$ and $v = (v_1, v_2)$ from $D(L)$ validity the following equality
$$(Lu, v)_{L^2} - (u, Lv)_{L^2} = (\gamma_1(u), \gamma_2(v))_H - (\gamma_2(u), \gamma_1(v))_H$$
can be easily verified. Now give any elements $f, g \in H$. Find the function $u = (u_1, u_2) \in D(L)$ such that
$$\gamma_1(u) = \frac{1}{i\sqrt{2}}(u_1(a) + u_2(b)) = f \text{ and } \gamma_2(u) = \frac{1}{\sqrt{2}}(u_1(a) - u_2(b)) = g$$
that is,
$$u_1(a) = (if + g)/\sqrt{2} \text{ and } u_2(b) = (if - g)/\sqrt{2}.$$
If choose these functions $u_1(t), u_2(t)$ in following form
$$u_1(t) = \int_{-\infty}^{t} e^{s-a} ds (if + g)/\sqrt{2}, \quad t < a,$$
$$u_2(t) = \int_{t}^{\infty} e^{b-t} ds (if - g)/\sqrt{2}, \quad t > b,$$
then it is clear that $(u_1, u_2) \in D(L)$ and $\gamma_1(u) = f$, $\gamma_2(u) = g$. □

Furthermore, using the method in [6] can be established the following result.

**Theorem 2.2:** If $\tilde{L}$ is a selfadjoint extension of the minimal operator $L_0$ in $L^2$, then it generates by the differential expression (2.1) and the boundary condition
$$u_2(b) = W u_1(a),$$
where $W: H \to H$ is a unitary operator. Moreover, the unitary operator $W$ in $H$ is determined uniquely by the extension $\tilde{L}$, i.e. $\tilde{L} = L(W)$ and vice versa.

## 3. The Spectrum of the Selfadjoint Extensions

In this section the structure of the spectrum of the selfadjoint extension $L_W$ in $L^2$ will be investigated.

First of all, we have to prove the following result.

**Theorem 3.1:** The point spectrum of selfadjoint extension $L_W$ is empty, i.e.
$$\sigma_p(L_W) = \emptyset.$$
**Proof:** Consider the following eigenvalue problem
$$\tilde{l}(u) = iu'(t) + \tilde{A}u(t) = \lambda u(t), \quad u \in L^2, \quad \lambda \in \mathbb{R}$$
$$u_2(b) = W u_1(a).$$
From this it is obtained that

$$u' = i(\tilde{A} - \lambda)u, \ u_2(b) = Wu_1(a), \ u \in L^2, \ \lambda \in \mathbb{R}.$$

The general solution of the last equation is
$$\begin{cases} u_1(t) = e^{i(\tilde{A}-\lambda)(t-a)}f, \ t < a \\ u_2(t) = e^{i(\tilde{A}-\lambda)(t-b)}g, \ t > b \\ u_2(b) = Wu_1(a), \ f, g \in H. \end{cases}$$

It is clear that for the $f \neq 0$, $g \neq 0$ the functions $u_1 \notin L^2(H, (-\infty, a))$, $u_2 \notin L^2(H, (b, \infty))$. So for every unitary operator $W$ we have $\sigma_p(L_W) = \emptyset$. □

Since residual spectrum of any selfadjoint operator in any Hilbert space is empty, then furthermore the continuous spectrum of selfadjoint extensions $L_W$ of the minimal operator $L_0$ is investigated.

Now it will be researched the resolvent of $L_W$ generated by the differential expression $\tilde{l}(\ )$ and the boundary condition
$$u_2(b) = Wu_1(a)$$
in the Hilbert space $L^2$, i.e.
$$\begin{cases} \tilde{l}(u) = iu'(t) + \tilde{A}u(t) = \lambda u(t) + f(t), \ u \in L^2, \ \lambda \in \mathbb{C}, \lambda_i = Im\lambda > 0 \\ u_2(b) = Wu_1(a) \end{cases} \quad (3.1)$$

Now we will be shown that the following function
$$u(\lambda; t) = (u_1(\lambda; t), u_2(\lambda; t)),$$
where
$$u_1(\lambda; t) = e^{-i(\lambda-\tilde{A})(t-a)}f_\lambda^* + i\int_t^a e^{-i(\lambda-\tilde{A})(t-s)}f(s)ds, \ t < a,$$
$$u_2(\lambda; t) = i\int_t^\infty e^{-i(\lambda-\tilde{A})(t-s)}f(s)ds, \ t > b,$$
$$f_\lambda^* = W^*\left(-i\int_b^\infty e^{-i(\lambda-\tilde{A})(b-s)}f(s)ds\right)$$

is a solution of the boundary value problems (3.1) in the Hilbert space $L^2$. For this, it is sufficient to show that
$$u_1(\lambda; t) \in L^2(H, (-\infty, a)),$$
$$u_2(\lambda; t) \in L^2(H, (b, +\infty))$$
for the case $\lambda_i > 0$. Indeed, in this case
$$\|f_\lambda^*\|_H^2 = \left\|-i\int_b^\infty e^{-i(\lambda-\tilde{A})(b-s)}f(s)ds\right\|_H^2 \leq \left(\int_b^\infty e^{\lambda_i(b-s)}\|f(s)\|_H ds\right)^2$$
$$\leq \left(\int_b^\infty e^{2\lambda_i(b-s)}ds\right)\left(\int_b^\infty \|f(s)\|_H^2 ds\right) = \frac{1}{2\lambda_i}\|f\|_{L^2(H,(b,+\infty))}^2 < \infty,$$

$$\left\|e^{-i(\lambda-\tilde{A})(t-a)}f_\lambda^*\right\|_{L^2(H,(-\infty,a))}^2 = \left\|e^{-i\lambda(t-a)}f_\lambda^*\right\|_{L^2(H,(-\infty,a))}^2 = \int_{-\infty}^a \left\|e^{-i\lambda(t-a)}f_\lambda^*\right\|_H^2 dt$$
$$= \int_{-\infty}^a e^{2\lambda_i(t-a)}dt \|f_\lambda^*\|_H^2 = \frac{1}{2\lambda_i}\|f_\lambda^*\|_H^2 < \infty$$

and
$$\left\|i\int_t^a e^{-i(\lambda-\tilde{A})(t-s)}f(s)ds\right\|_{L^2(H,(-\infty,a))}^2 \leq \int_{-\infty}^a \left(\int_t^a e^{\lambda_i(t-s)}\|f(s)\|_H ds\right)^2 dt$$

$$\leq \int_{-\infty}^{a}\left(\int_{t}^{a} e^{\lambda_i(t-s)}ds\right)\left(\int_{t}^{a} e^{\lambda_i(t-s)}\|f(s)\|^2 ds\right)dt$$

$$= \frac{1}{\lambda_i}\int_{-\infty}^{a}\int_{t}^{a} e^{\lambda_i(t-s)}\|f(s)\|^2 ds dt = \frac{1}{\lambda_i}\int_{-\infty}^{a}\left(\int_{-\infty}^{s} e^{\lambda_i(t-s)}\|f(s)\|^2 dt\right)ds$$

$$= \frac{1}{\lambda_i}\int_{-\infty}^{a}\left(\int_{-\infty}^{s} e^{\lambda_i(t-s)}dt\right)\|f(s)\|^2 ds = \frac{1}{\lambda_i^2}\int_{-\infty}^{a}\|f(s)\|^2 ds$$

$$= \frac{1}{\lambda_i^2}\|f\|_{L^2(H,(-\infty,a))}^2 < \infty.$$

Furthermore,

$$\left\|i\int_{t}^{\infty} e^{-i(\lambda-\tilde{A})(t-s)}f(s)ds\right\|_{L^2(H,(b,+\infty))} \leq \int_{b}^{\infty}\left(\int_{t}^{\infty} e^{\lambda_i(t-s)}\|f(s)\|_H ds\right)^2 dt$$

$$\leq \int_{b}^{\infty}\left(\int_{t}^{\infty} e^{\lambda_i(t-s)}ds\right)\left(\int_{t}^{\infty} e^{\lambda_i(t-s)}\|f(s)\|^2 ds\right)dt$$

$$= \frac{1}{\lambda_i}\int_{b}^{\infty}\left(\int_{t}^{\infty} e^{\lambda_i(t-s)}\|f(s)\|^2 ds\right)dt = \frac{1}{\lambda_i}\int_{b}^{\infty}\left(\int_{b}^{s} e^{\lambda_i(t-s)}\|f(s)\|^2 dt\right)ds$$

$$= \frac{1}{\lambda_i}\int_{b}^{\infty}\left(\int_{b}^{s} e^{\lambda_i(t-s)}dt\right)\|f(s)\|^2 ds = \frac{1}{\lambda_i^2}\left(\int_{b}^{\infty}(1-e^{\lambda_i(b-s)})\|f(s)\|^2 ds\right)$$

$$\leq \frac{1}{\lambda_i^2}\|f\|_{L^2(H,(b,+\infty))}^2 < \infty.$$

From above calculations imply that $u_1(\lambda;t) \in L^2(H,(-\infty,a))$, $u_2(\lambda;t) \in L^2(H,(b,+\infty))$ for $\lambda \in \mathbb{C}$, $\lambda_i = Im\lambda > 0$. On the other hand it can be to easy to verify that $u(\lambda;t)$ is a solution of the boundary value problem (3.1).

In the case when $\lambda \in \mathbb{C}$, $\lambda_i = Im\lambda < 0$ solution of the boundary value problem
$$L_W u = iu' + \tilde{A}u = \lambda u + f, \quad u = (u_1, u_2), \quad f \in L^2$$
$$u_2(b) = Wu_1(a),$$
where $W$ is a unitary operator in $H$, is in the form $u(\lambda;t) = (u_1(\lambda;t), u_2(\lambda;t))$,

$$\begin{cases} u_1(\lambda;t) = -i\int_{-\infty}^{t} e^{-i(\lambda-\tilde{A})(t-s)}f(s)ds, & t < a \\ u_2(\lambda;t) = e^{-i(\lambda-\tilde{A})(t-b)}g_\lambda^* - i\int_{b}^{t} e^{-i(\lambda-\tilde{A})(t-s)}f(s)ds, & t > b, \end{cases}$$

$$g_\lambda^* = W\left(-i\int_{-\infty}^{a} e^{-i(\lambda-\tilde{A})(a-s)}f(s)ds\right).$$

In first prove that $u(\lambda;t) \in L^2$. In this case

$$\|u_1(\lambda;t)\|_{L^2(H,(-\infty,a))}^2 = \int_{-\infty}^{a}\left\|-i\int_{-\infty}^{t} e^{-i(\lambda-\tilde{A})(t-s)}f(s)ds\right\|_H^2 dt$$

$$\leq \int_{-\infty}^{a}\left(\int_{-\infty}^{t} e^{\lambda_i(t-s)}ds\right)\left(\int_{-\infty}^{t} e^{\lambda_i(t-s)}\|f(s)\|_H^2 ds\right)dt$$

$$= \frac{1}{|\lambda_i|}\int_{-\infty}^{a}\int_{-\infty}^{t} e^{\lambda_i(t-s)}\|f(s)\|_H^2\, ds\, dt$$

$$= \frac{1}{|\lambda_i|}\int_{-\infty}^{a}\left(\int_{s}^{a} e^{\lambda_i(t-s)}\|f(s)\|_H^2 dt\right)ds = \frac{1}{|\lambda_i|}\int_{-\infty}^{a}\left(\int_{s}^{a} e^{\lambda_i(t-s)}dt\right)\|f(s)\|_H^2 ds$$

$$= \frac{1}{|\lambda_i|^2}\int_{-\infty}^{a}(1-e^{\lambda_i(a-s)})\|f(s)\|_H^2\, ds \leq \frac{1}{|\lambda_i|^2}\|f\|_{L^2(H,(-\infty,a))}^2 < \infty,$$

$$\|g_\lambda^*\|_H^2 = \left\|-i\int_{-\infty}^{a} e^{-i(\lambda-\tilde{A})(a-s)}f(s)ds\right\|_H^2 \leq \left(\int_{-\infty}^{a} e^{\lambda_i(a-s)}\|f(s)\|_H ds\right)^2$$

$$\leq \left(\int_{-\infty}^{a} e^{2\lambda_i(a-s)}ds\right)\left(\int_{-\infty}^{a}\|f(s)\|_H^2\, ds\right)$$

$$= \frac{1}{2|\lambda_i|}\|f\|_{L^2(H,(-\infty,a))}^2 < \infty,$$

$$\left\|e^{-i(\lambda-\tilde{A})(t-b)}g_\lambda^*\right\|_{L^2(H,(b,+\infty))}^2 \leq \int_{b}^{\infty} e^{2\lambda_i(t-b)}dt\, \|g_\lambda^*\|_H^2 = \frac{1}{2|\lambda_i|}\|g_\lambda^*\|_H^2$$

$$\leq \frac{1}{4|\lambda_i|^2}\|f\|_{L^2(H,(b,+\infty))}^2 < \infty,$$

and

$$\left\|-i\int_{b}^{t} e^{-i(\lambda-\tilde{A})(t-s)}f(s)ds\right\|_{L^2(H,(b,+\infty))}^2 \leq \int_{b}^{\infty}\left(\int_{b}^{t} e^{\lambda_i(t-s)}\|f(s)\|_H ds\right)^2 dt$$

$$\leq \int_{b}^{\infty}\left(\int_{b}^{t} e^{\lambda_i(t-s)}ds\right)\left(\left(\int_{b}^{t} e^{\lambda_i(t-s)}\|f(s)\|_H^2 ds\right)\right)dt$$

$$= \int_{b}^{\infty}\left(\frac{1}{\lambda_i}(1-e^{\lambda_i(t-b)})\right)\left(\int_{b}^{t} e^{\lambda_i(t-s)}\|f(s)\|_H^2 ds\right)dt$$

$$\leq \frac{1}{|\lambda_i|}\int_{b}^{\infty}\left(\int_{b}^{t} e^{\lambda_i(t-b)}\|f(s)\|_H^2 ds\right)dt$$

$$= \frac{1}{|\lambda_i|}\int_{b}^{\infty}\left(\int_{s}^{\infty} e^{\lambda_i(t-s)}\|f(s)\|_H^2 dt\right)ds$$

$$= \frac{1}{|\lambda_i|}\int_{b}^{\infty}\left(\int_{s}^{b} e^{\lambda_i(t-s)}dt\right)\|f(s)\|_H^2 ds$$

$$= \frac{1}{|\lambda_i|^2}\|f\|_{L^2(H,(b,+\infty))}^2 < \infty.$$

The above simple calculations are shown that $u_1(\lambda;\cdot) \in L^2(H,(-\infty,a))$, $u_2(\lambda;\cdot) \in L^2(H,(b,+\infty))$, i.e. $u(\lambda;\cdot) \in L^2$ in the case when $\lambda \in \mathbb{C}$, $\lambda_i = Im\lambda < 0$.

On the other hand it can be verified that the function $u(\lambda;\cdot)$ satisfy the equation $\tilde{l}u(\lambda;\cdot) = iu'(\lambda;\cdot) + \tilde{A}u(\lambda;\cdot) = \lambda u(\lambda;\cdot) + f$ and $u_2(b) = Wu_1(a)$.

Hence the following result has been proved.

**Theorem 3.2:** For the resolvent set $\rho(L_W)$ is valid
$$\rho(L_W) \supset \{\lambda \in \mathbb{C}: Im\lambda \neq 0\}.$$

Now will be researched continuous spectrum $\sigma_c(L_W)$ of the extension $L_W$. For the $\lambda \in \mathbb{C}, \lambda_i = Im\lambda > 0$ the norm of resolvent operator $R_\lambda(L_W)$ of the $L_W$ is in form

$$\|R_\lambda(L_W)f(t)\|_{L^2}^2 = \left\| e^{-i(\lambda-\tilde{A})(t-a)}f_\lambda^* + i\int_t^a e^{-i(\lambda-\tilde{A})(t-s)}f_1(s)ds \right\|_{L^2(H,(-\infty,a))}^2$$

$$+ \left\| i\int_t^\infty e^{-i(\lambda-\tilde{A})(t-s)}f_2(s)ds \right\|_{L^2(H,(b,+\infty))}^2$$

$f \in L^2$, $f = (f_1, f_2)$.

Then it is clear that for any $f = (f_1, f_2) \in L^2$ is true
$$\|R_\lambda(L_W)f(t)\|_{L^2}^2 \geq \left\| i\int_t^\infty e^{-i(\lambda-\tilde{A})(t-s)}f_2(s)ds \right\|_{L^2(H,(b,+\infty))}^2.$$

The vector functions $f^*(\lambda;t)$ in form $f^*(\lambda;t) = (0, e^{-i(\bar{\lambda}-\tilde{A})t}f)$, $\lambda \in \mathbb{C}$, $\lambda_i = Im\lambda > 0$, $f \in H$ belong to $L^2$. Indeed,
$$\|f^*(\lambda;t)\|_{L^2}^2 = \int_b^\infty \|e^{-i(\bar{\lambda}-\tilde{A})t}f\|_H^2 dt = \int_b^\infty e^{-2\lambda_i t}dt \|f\|_H^2 = \frac{1}{2\lambda_i}e^{-2\lambda_i b} < \infty.$$

For the such functions $f^*(\lambda;\cdot)$ we have

$$\|R_\lambda(L_W)f^*(\lambda;t)\|_{L^2(H,(b,+\infty))}^2 \geq \left\| i\int_t^\infty e^{-i(\lambda-\tilde{A})(t-s)}e^{-i(\bar{\lambda}-\tilde{A})s}fds \right\|_{L^2(H,(b,+\infty))}^2$$

$$= \left\| \int_t^\infty e^{-i\lambda t}e^{-2\lambda_i s}e^{i\tilde{A}t}fds \right\|_{L^2(H,(b,+\infty))}^2 = \left\| e^{-i\lambda t}e^{i\tilde{A}t}\int_t^\infty e^{-2\lambda_i s}fds \right\|_{L^2(H,(b,+\infty))}^2$$

$$= \left\| e^{-i\lambda t}\int_t^\infty e^{-2\lambda_i s}ds \right\|_{L^2(H,(b,+\infty))}^2 \|f\|_H^2 = \frac{1}{4\lambda_i^2}\int_b^\infty e^{-2\lambda_i t}dt \|f\|_H^2 = \frac{1}{8\lambda_i^3}e^{-2\lambda_i b}\|f\|_H^2.$$

From this
$$\|R_\lambda(L_W)f^*(\lambda;\cdot)\|_{L^2} \geq \frac{e^{-\lambda_i b}}{2\sqrt{2}\lambda_i\sqrt{\lambda_i}}\|f\|_H = \frac{1}{2\lambda_i}\|f^*(\lambda;\cdot)\|_{L^2}$$

i.e. for $\lambda_i = Im\lambda > 0$ and $f \neq 0$ are valid
$$\frac{\|R_\lambda(L_W)f^*(\lambda;\cdot)\|_{L^2}}{\|f^*(\lambda;\cdot)\|_{L^2}} \geq \frac{1}{2\lambda_i}.$$

On the other hand it is clear that
$$\|R_\lambda(L_W)\| \geq \frac{\|R_\lambda(L_W)f^*(\lambda;\cdot)\|_{L^2}}{\|f^*(\lambda;\cdot)\|_{L^2}}, f \neq 0.$$

Consequently, we have
$$\|R_\lambda(L_W)\| \geq \frac{1}{2\lambda_i} \text{ for } \lambda \in \mathbb{C}, \lambda_i = Im\lambda > 0.$$

Actually, this has been proved the following claim.

**Theorem 3.3:** Continuous spectrum of the extension $L_W$ in form
$$\sigma_c(L_W) = \mathbb{R}.$$

**Example:** By the last theorem the spectrum of following boundary value problem
$$i\frac{\partial u(t,x)}{\partial t} - \frac{\partial^2 u(t,x)}{\partial x^2} = f(t,x),\ |t| > 1,\ x \in [0,1],$$
$$u(1,x) = e^{i\varphi}u(-1,x),\ \varphi \in [0,2\pi),$$
$$u'_x(t,0) = u'_x(t,1) = 0,\ |t| > 1,$$
in the space $L^2\big((-\infty,-1) \times (0,1)\big) \oplus L^2\big((1,\infty) \times (0,1)\big)$ is continuous and coincides with $\mathbb{R}$.

Later on, note that another approach has been given in [9] for the singular differential operators for n-th order in scalar case.